\documentclass[12 pt]{article}
\usepackage{amsfonts,amsmath,color,amssymb, amsthm}
\usepackage{latexsym}
\title{Covering $n$-Permutations with $(n+1)$-Permutations}
\author{Taylor F.~Allison\\
Department of Mathematics\\
North Carolina State University\\
Raleigh, NC, USA \and
Anant P.~Godbole\\
Department of Mathematics \& Statistics\\
East Tennessee State University\\
Johnson City, TN, USA\and
Kathryn M.~Hawley\\
Department of Mathematics\\
Harvey Mudd College\\
Claremont, CA, USA\and
Bill Kay\\
Department of Mathematics\\
University of South Carolina,
Columbia, SC, USA.}
\begin{document}
\def\n{\emptyset}
\def\lr{\left(}
\def\rr{\right)}
\def\lc{\left\{}
\def\rc{\right\}}
\def\vp{\varphi}
\def\ca{\mathcal A}
\def\e{\mathbb E}
\def\p{\mathbb P}
\def\v{\mathbb V}
\newcommand{\marginal}[1]{\marginpar{\raggedright\scriptsize #1}}
\newcommand{\BK}{\marginal{BK}}
\newcommand{\AG}{\marginal{AG}}
\newcommand{\KH}{\marginal{KH}}
\newcommand{\TA}{\marginal{TA}}
\newtheorem{thm}{Theorem}
\newtheorem{lm}[thm]{Lemma}
\newtheorem{rem}[thm]{Remark}
\newtheorem{exam}[thm]{Example}
\newtheorem{prop}[thm]{Proposition}
\newtheorem{defn}[thm]{Definition}
\newtheorem{cm}[thm]{Claim}
\newtheorem{conj}[thm]{Conjecture}
\def\ep{\varepsilon}
\def\tv{{d_{{\rm TV}}}}
\def\P{{\rm Po}}
\def\S{S}
\def\k{\kappa}
\def\cl{{\cal L}}
\def\ck{{\cal K}}
\def\p{\mathbb P}
\def\v{\mathbb V}
\def\e{\mathbb E}
\def\z{\mathbb Z}
\def\l{\lambda}
\maketitle
\begin{abstract} Let $S_n$ be the set of all permutations on $[n]:=\{1,2,\ldots,n\}$. We denote by $\k_n$ the smallest cardinality of a subset $\ca$ of $\S_{n+1}$ that ``covers" $\S_n$, in the sense that each $\pi\in\S_n$ may be found as an order-isomorphic subsequence of some $\pi'$ in $\ca$.  What are general upper bounds on $\k_n$?  If we randomly select $\nu_n$ elements of $\S_{n+1}$, when does the probability that they cover $\S_n$ transition from 0 to 1?  Can we provide a fine-magnification analysis that provides the ``probability of coverage"  when $\nu_n$ is around the level given by the phase transition?   In this paper we answer these questions and raise others.
\end{abstract}
\section{Introduction}  The problem discussed in this paper was posed by Prof.~Robert Brignall during the Open Problem Session at the {\it International Permutation Patterns Conference} held at California State Polytechnic University in June 2011.  The conference webpage may be found at 

\centerline{\tt http://www.calpoly.edu/$\sim$math/PP2011/}.  

\noindent Brignall asked, ``If $S_n$ denotes the set of all permutations on $[n]:=\{1,2,\ldots,n\}$, what can we say about $\k_n$, the smallest cardinality of a subset $\ca$ of $\S_{n+1}$ that ``covers" $\S_n$, in the sense that each $\pi\in\S_n$ may be found as an order-isomorphic subsequence of some $\pi'$ in $\ca$."  Specifically, he asked about bounds, exact values, and asymptotics.  Small values are easy to calculate; e.g., it is easy to see that $\k_1=1, \k_2=1,$ and the permutation set $\{1342, 4213\}$ reveals that $\k_3=2$ -- but the situation rapidly gets out of precise control.  

In Section 2, we provide a trivial pigeonhole lower bound on $\k_n$ and then use the ``method of alterations" \cite{as} to derive an upper general bound on $\k_n$ that contains a logarithmic factor that often shows up when using such arguments in covering-type situations; see, e.g. the general upper bound on the size of covering designs that was proved by Erd\H os and Spencer in their early work \cite{es}.  We continue by showing that the {\it second and subsequent} coverings of the $n$-permutations are accomplished in linear $\log\log$ time, in  a result that is reminiscent of the ones in \cite{gtv} (covering designs) and \cite{gss} ($t$-covering arrays).  Lastly, in Section 2, we make comparisons to the development in Spencer \cite{sp} to produce evidence that the asymptotic value of $\k_n$ is (upto a $O(1)$ or perhaps $1+o(1)$ factor) the same as that given by the lower bound; we conjecture that this is true.

In Section 3, we switch to a different approach, asking the question ``If we randomly select $\nu_n$ elements of $\S_{n+1}$, when does the probability that they cover $\S_n$ transition from asymptotically zero 0 to asymptotically 1?"  Our main result in this area, Theorem 6, is proved using the Janson exponential inequality (\cite{as}).  The threshold in Theorem 6 contains a small gap, but the result  is the best possible, as we show in Section 4, where the Stein-Chen method of Poisson approximation is used to prove, that in the regime of $\nu_n$s that are ``in the gap," the number of uncovered $n$-permutations $X$ has a Poisson distribution with finite mean, and thus $\p(X=0)$ is a  finite constant that is bounded away from zero and one. 
\section{Bounds}
Our first preliminary result provides a formula for the number $c(n,\pi)$ of permutations in $\S_{n+1}$ that cover a fixed $\pi\in \S_n$.
\begin{lm}  Let $c(n,\pi)$ denote the number of permutations in $\S_{n+1}$ that cover a fixed $\pi\in \S_n$.  Then $c(n,\pi)=c(n,\pi')=n^2+1$ for each $\pi, \pi'\in \S_n$.
\end{lm}
\begin{proof}
It is clear that any permutation pattern $\pi\in\S_n$ may be realized in ${{n+1}\choose{n}}=n+1$ ways, one for any choice of $n$ numbers from $\{1,2,\ldots,n+1\}$.  Arrange these ways lexicographically (for example if $n=3$, we can realize the pattern 132 as 132, 142, 243, and 143, or, lexicographically, as 132, 142, 143, 243).  Note that the $r$th and $(r+1)$st lex-orderings of $\pi$ differ in a single bit.  Now, given any realization of $\pi$, the $(n+1)$st letter may clearly be inserted in $(n+1)$ ways to create an $(n+1)$-covering permutation; however, for any $1\le r\le n-1$, the list of covering $(n+1)$-permutations for the $r$th and $(r+1)$st lex-orderings have an overlap of magnitude 2, corresponding to whether the $(n+1)$st letter is inserted before or after the non-matching bit.  Thus $c(n,\pi)=c(n,\pi')=(n+1)^2-2n=n^2+1$, as asserted.
\end{proof}
Lemma 1 can now be used to prove
\begin{thm}
\[\frac{(n+1)!}{n^2}(1+o(1))\le\k_n\le\frac{\log n}{n^2}(n+1)!(1+o(1)).\]
\end{thm}
\begin{proof}
The lower bound is elementary.  Each of the $\k_n$ permutations in a covering, ``takes care," with repetition, of ${{n+1}\choose{n}}=n+1$ $n$-permutations.  Since we have a covering, clearly, $\k_n(n+1)\ge n!$, or
$$
\k_n\ge\frac{n!}{(n+1)}=\frac{(n+1)!}{n^2}\lr1-\frac{2n+1}{(n+1)^2}\rr.
$$
For the upper bound, we use the method of alterations \cite{as} as follows:  Choose a random number $Y$ of $(n+1)$-permutations by ``without replacement" sampling.  The expected number of uncovered $n$-permutations, by Lemma 1 and linearity of expectation, is $$\e(X)=n!\frac{{{(n+1)!-n^2-1}\choose{Y}}}{{{(n+1)!}\choose{Y}}}.$$ We choose a realization with $X=X_Y\le\e(X)$ and cover these with at most 
\[n!\lr1-\frac{n^2+1}{(n+1)!}\rr^Y\le n!\exp\{-Y(n^2+1)/(n+1)!\}\]
additional $(n+1)$-permutations, yielding, for any initial size $Y$, a covering with at most 
\[Y+n!\exp\{-Y(n^2+1)/(n+1)!\}\]
members.  Minimizing over $Y$ yields an initial choice of size 
\[\frac{(n+1)!}{(n^2+1)}\log\lr\frac{n^2+1}{n+1}\rr,\]
and an upper bound of
$$
\k_n\le\frac{(n+1)!}{(n^2+1)}\lr1+\log\lr\frac{n^2+1}{n+1}\rr\rr=\frac{\log n}{n^2}(n+1)!(1+o(1)),
$$
as claimed.
\end{proof}
We next ask how many more $(n+1)$-permutations are required to cover the $n$ permutations multiple times.  Here the situation is often nuanced, and leads to the question as to whether the logarithmic factor, present in the upper bound for the first covering, is extraneous.  Typically, in other covering contexts, we find that the second and subsequent coverings need an appropriately normalized $\log\log n$ additional elements in the cover; see, e.g. \cite{gtv} for a {\it covering design} analogy, \cite {gss} for an occurrence of this phenomenon in {\it $t$-covering arrays}, and \cite{fs},  \cite {mw} for the $\log\log$ behavior in the {\it coupon collection problem}.  We briefly describe the parallel in the context of covering designs:  A collection $\ca$ of sets of size $k$ of $[n]$ is said to form a $t$-covering design if each $t$-set is contained in at least one $k$-set in $\ca$.  If $m(n,k,t)$ denotes the smallest size of a $t$-covering design $\ca$ then it is clear that $m(n,k,t)\ge{n\choose t}/{k\choose t}$; Erd\H os and Spencer proved in \cite{es} that $\forall n,k,t$,
\[m(n,k,t)\le\frac{{n\choose t}}{{k\choose t}}\lr1+\log {k\choose t}\rr;\]
it was shown furthermore in \cite{gtv} that the minimum number $m(n,k,t,\l)$ of $k$-sets needed to cover each $t$-set $\l$ times satisfied
\[m(n,k,t,\l)\le\frac{{n\choose t}}{{k\choose t}}\lr1+\log {k\choose t}+(\l-1)\log\log{k\choose t}+O(1)\rr,\]
$n,k,t\to\infty$.  This was the $\log\log$ result.  Also, the Erd\H os-Hanani conjecture, namely that for fixed $k,t$,
\[\lim_{n\to\infty}\frac{m(n,k,t)}{{n\choose t}}=\frac{1}{{k\choose t}}\]
was proved by R\"odl \cite{ro} and, later, by Spencer \cite{sp}.  This showed that the logarithmic factor in the Erd\H os-Spencer bound could be asymptotically dispensed with.  Finally, see \cite{gj} for a corresponding threshold result.  It is these questions we seek to address, in our context, in the rest of this section and the next.
\begin{thm} Let $\k_{n,\l}$ denote the minimum number of $(n+1)$-permutations needed to cover each $n$-permutation $\l\ge2$ times.  Then,
\[\k_{n,\l}\le\frac{(n+1)!}{n^2}\lr\log n+(\l-1)\log\log n+O(1)\rr.\]
\end{thm}
\begin{proof}
We first choose an unspecified number $Y$ of $(n+1)$-permutations randomly and {\it with} replacement.  This might lead to replication with very small probability, but the proof is far more streamlined - and can easily be adapted to the case where we choose $Y$ distinct permutations.  The probability that any permutation is covered just $j$ times; $0\le j\le\l-1$ is
\[{Y\choose j}\lr\frac{n^2+1}{(n+1)!}\rr^j\lr1-\frac{n^2+1}{(n+1)!}\rr^{Y-j},\]
so that the expected number of such permutations is
\[\e(X_j)=n!{Y\choose j}\lr\frac{n^2+1}{(n+1)!}\rr^j\lr1-\frac{n^2+1}{(n+1)!}\rr^{Y-j}; 0\le j\le\l-1.\]
We can cover each such permutation in any {\it ad hoc} way by choosing $\l-j$ additional $(n+1)$-permutations, and so the number $Z_Y$ of $(n+1)$-permutations in this successful $\l$-covering is
\[Z_Y=Y+\sum_{j=0}^{\l-1}(\l-j)X_j,\]
and thus
\[\e(Z_Y)=Y+n!\sum_{j=0}^{\l-1}(\l-j){Y\choose j}\lr\frac{n^2+1}{(n+1)!}\rr^j\lr1-\frac{n^2+1}{(n+1)!}\rr^{Y-j}\]
for any initial choice of $Y$ permutations. Set 
\[p=\frac{n^2+1}{(n+1)!}; q=1-p,\] and, given $Y$ Bernoulli trials with success probability $p$, denote the cumulative and point binomial probabilities by $B(Y,p,k)$ and $b(Y,p,k)$ respectively, $0\le k\le Y$.  It is easy to verify that for $Yp\ge\l$
\begin{eqnarray}\e(Z_Y)&=&Y+n!\l B(Y,p,\l-1)-n!YpB(Y-1,p,\l-2)\nonumber\\
&=&Y+n!\l b(Y,p,\l-1)+n!\lc\l B(Y,p,\l-2)-YpB(Y-1,p,\l-2)\rc\nonumber\\
&\le&Y+n!\l b(Y,p,\l-1)+n!\l\{B(Y,p,\l-2)-B(Y-1,p,\l-2)\}\nonumber\\
&=&Y+n!\l b(Y,p,\l-1)+n!\l\sum_{j=0}^{\l-2}{{Y-1}\choose{j}}p^jq^{Y-1-j}\lr\frac{Yq}{Y-j}-1\rr\nonumber\\
&\le&Y+n!\l b(Y,p,\l-1)\nonumber\\
&\le&Y+n!\l\lr\frac{Yp}{q}\rr^{\l-1}\frac{e^{-pY}}{(\l-1)!}.\end{eqnarray}
We do not attempt to optimize carefully in (1); rather we set
\[Y=\frac{(n+1)!}{n^2+1}(\log n+(\l-1)\log\log n)\] and see that the second term $T_2$ in (1) reduces as
\begin{eqnarray*}
T_2&=&n!\l\lr\frac{\log n+(\l-1)\log\log n}{q}\rr^{\l-1}\frac{\exp\{-\log n-(\l-1)\log\log n\}}{(\l-1)!}\\
&=&\frac{(n+1)!}{n^2+1}\frac{\l}{(\l-1)!}\frac{1}{(\log n)^{\l-1}}\lr\frac{\log n+(\l-1)\log\log n}{q}\rr^{\l-1}(1+o(1))\\
&=&\frac{(n+1)!}{n^2+1}\frac{\l}{(\l-1)!}(1+o(1)),\end{eqnarray*}
so that (1) yields
\[\e(Z_Y)\le \frac{(n+1)!}{n^2+1}\lr\log n+(\l-1)\log\log n+(1+o(1))\frac{\l}{(\l-1)!}\rr.\]
Finally, we choose a sample outcome for which $\vert \ca\vert\le \e(Z_Y)$ to complete the proof.
\end{proof}
We now describe the hypergraph formulation of \cite{ps} that was used in  \cite{sp} to prove the Erd\H os-Hanani conjecture using a method that involved branching processes, dynamical algorithms, hypergraph theory, and differential equations.  In this formulation the vertices of the hypergraph consisted of the ensemble of $t$-sets; for us they would be the class of permutations in $S_n$.  The edges in \cite{sp} were the collections of $t$-subsets of the $k$-sets, so that the hypergraph was ${k\choose t}$ uniform.  If analogously, we let edges be the set of $n$-permutations covered by an $(n+1)$-permutation, then the hypergraph is no longer uniform.  It is not too hard to prove, however, that each $(n+1)$-permutation $\pi$ covers $n+1-s_\pi$ $n$-permutations, where $s_\pi$ is the number of {\it successions} in $\pi$, where a succession is defined as an episode $\pi(i+1)=\pi(i)\pm1$.  Moreover, we know \cite{gs} that the number of successions in a random permutation is approximately Poisson with parameter $\sim2$, so that it is reasonable to assert that most hypergraph edges consist of $n-O(1)$ vertices.  This is the first deviation from the Pippenger model, which we consider to be {\it not too serious} insofar as the lack of uniformity of the hypergraph is concerned but {\it rather serious} due to the fact that the uniformity level $n-O(1)$ is not finite. Lemma 1 above shows that the degree of each vertex is $O(n^2)$, and we will prove in Lemma 5 below that the codegree of two vertices $\pi$ and $\pi'$ is at most $O(1)$, so that the codegree is an order of magnitude smaller than the degree.  This is good.  The above problems with the hypergraph formulation notwithstanding, we make the following conjecture:
\begin{conj}  For some constant $A$,
\[\lim\sup_{n\to\infty}\frac{\k_n}{(n+1)!/n^2}=A,\]
and possibly $A\le2$.
\end{conj}
\section{Thresholds}
If each permutation in $S_{n+1}$ is randomly and independently picked with probability $p$, we will show that the probabilistic zero-one threshold for coverage of $S_n$ is at the level $p=\log n/n$, i.e., at $n$ times the level given by the upper bound in Theorem 2.  Towards this end we prove the following result, of interest in its own right.
\begin{lm}  For any $\pi\in\S_n$, the set
\[{\cal J}_\pi:=\{\pi'\in\S_n: \pi\ {\rm and}\ \pi'\ {\rm can\ be\ jointly\ covered\ by\ }\rho\in\S_{n+1}\}\]
has cardinality at most $n^3$.  Moreover, for any $\pi, \pi'\in\S_{n}$, the cardinality of 
\[{\cal C}_{\pi, \pi'}:=\{\rho\in\S_{n+1}:\rho\ {\rm covers\ both\ } \pi\ {\rm and\ }\pi'\}\]
is at most $4$.
\end{lm}
\begin{proof} 
Fix $\pi$.  For an $(n+1)$-permutation to be able to successfully cover another $\pi'\in S_n$ (in addition to $\pi$), $\pi$ must contain an $(n-1)$-subpattern of $\pi'$.  This subpattern may be present in ${{n}\choose{n-1}}=n$ possible positions of $\pi$, and can be represented, using the numbers $\{1,2,\ldots,n\}$, in $n$ ways.  Finally, the $n$th letter of $\pi'$ can be inserted into this subpattern in $n$ ways.  This proves the first part of the lemma; we observe that the $n^3$ bound is probably not best and try much harder   to get the correct bound in the second part.

To prove the second part, we note that $\rho\in\S_{n+1}$ is able to successfully cover two co-coverable permutations $\pi$ and $\pi'$ only if there exists an insertion of the symbol $\emptyset$ into the sequences $\pi=(\pi_1,\pi_2,\ldots,\pi_n)$ and $\pi'=(\pi'_1,\pi'_2,\ldots,\pi'_n)$, so that when the augmented sequences $\pi_{\emptyset}=\sigma=(\sigma_1,\sigma_2,\ldots,\sigma_{n+1})$ and $ \pi'_{\n}=\sigma'=(\sigma'_1,\sigma'_2,\ldots,\sigma'_{n+1})$ are laid atop each other, the following conditions hold:

\noindent (i) There is a minimum index $h$ such that $$\sigma_{i_{n+1}}=\sigma'_{i_{n+1}}=n, \sigma_{i_{n}}=\sigma'_{i_{n}}=n-1,\ldots, \sigma_{i_h}=\sigma'_{i_h}=h-1,$$
(we set $h=n+2$ if no such index exists);

\noindent (ii) Either $\ell=0$ or else there is a maximum index $\ell$ such that 
$$\sigma_{i_1}=\sigma'_{i_1}=1, \sigma_{i_{2}}=\sigma'_{i_{2}}=2,\ldots, \sigma_{i_\ell}=\sigma'_{i_\ell}=\ell;$$
(iii) Without loss of generality (this is the ``$\sigma\le\sigma'$" solution),
$$\sigma_{i_{\ell+1}}=\emptyset; \sigma'_{i_{\ell+1}}=\ell+1;$$
$$\sigma_{i_{h-1}}=h-2; \sigma'_{i_{h-1}}=\emptyset;$$
(iv) and, lastly, $$\sigma_{i_j}=j-1;\sigma'_{i_j}=j; \ell+2\le j\le h-2.$$
The basic reason for this is that in the covering $(n+1)$-permutation $\rho$, the numbers 1 and $n+1$ may only play the role of 1 and $n$ respectively in either of the permutations $\pi$, $\pi'$, while any other number $t$ in $\rho$ may play the role of either $t$ or $t-1$. Let us give an example.  If $\pi=126983745; \pi'=127938645$, then a legal insertion of the $\emptyset$ symbols, with $h=10$ and $l=5$,  is as follows:
\begin{eqnarray*}
\sigma:&&\quad1\quad2\quad6\quad9\quad8\quad3\quad7\quad\emptyset\quad4\quad5\\
\sigma':&&\quad1\quad2\quad7\quad9\quad\emptyset\quad3\quad8\quad6\quad4\quad5\\
\rho:&&\quad1\quad2\quad7\quad10\quad9\quad3\quad8\quad6\quad4\quad5,\\
\end{eqnarray*}
whereas a ``bad" insertion would be, e.g., 
\begin{eqnarray*}
\sigma:&&\quad1\quad2\quad6\quad9\quad\emptyset\quad8\quad3\quad7\quad4\quad5\\
\sigma':&&\quad1\quad2\quad7\quad9\quad3\quad8\quad\emptyset\quad6\quad4\quad5,\\
\end{eqnarray*}
 since, starting with 10 representing the two 9's and 9 representing the two 8's -- we find that there is no position, in $\rho$, where the 8 can be placed.  

Notice also that if the arrangements $\sigma$ and $\sigma'$ are legal, then they may be laid atop each other, by simultaneous rearrangement, so that each is left-right non-decreasing; and the covering $(n+1)$-permutation is $id_{n+1}$.  If the arrangements are not legal, then such a rearrangement leads to sequences that are not both non-decreasing.  For instance, the legal example above can be rearranged as follows:

\begin{eqnarray*}
\sigma:&&\quad1\quad2\quad3\quad4\quad5\quad\emptyset\quad6\quad7\quad8\quad9\\
\sigma':&&\quad1\quad2\quad3\quad4\quad5\quad6\quad7\quad8\quad\emptyset\quad9\\
\end{eqnarray*}

We shall prove that for any two co-coverable permutations $\pi,\pi'$, there are at most two legal insertions of the $\emptyset$ symbols, for each of which there can be at most two covering permutations.  First, align any two co-coverable $\pi$ and $\pi'$ in $n-1$ spots, and introduce two $\emptyset$ symbols so that (i) through (iv) above are satisfied.
It is clear that this can be done, because
there exists a $\rho$ which covers both $\pi$ and $\pi'$, and deletion of two symbols from $\rho$ yields a $(n-1)$-subpattern shared by $\pi, \pi'$. The question is when this process is non-unique.  
% A study of all cases reveals that the only case where moving one or both $\emptyset$ symbols {\it does not} 
% lead to one of more of the conditions (i) through (iv) being violated is when the alignment contains the sub-alignment
%\begin{eqnarray*}
%\sigma:&&\quad a\quad b \quad\emptyset\\
%\sigma':&&\quad\emptyset\quad b\quad a,\\
%\end{eqnarray*}
%which can be rearranged as
%\begin{eqnarray*}
%\sigma:&&\quad \emptyset \quad a \quad b\\
%\sigma':&&\quad b\quad a\quad\emptyset,\\
%\end{eqnarray*}
%or vice versa.  In this case the corresponding elements of a covering $\rho$ can be (1) $a+1, b+1, a$; (2) $a,b+1,a+1$; (3) $b,a+1,b+1$; or (4) $b+1,a+1,b$.  This completes the proof.
It will be useful to introduce the notion of {\it longest common match}, or LCM, between two permutations.  Given two co-coverable permutations it is easy to determine which one is the ``larger" or ``dominant" one.  Call any dominant permutation $\pi'$ and the other $\pi$ ($\pi'$ may not be unique, as in Example (ii) below.)  When $\pi$ is laid atop $\pi'$, we identify an LCM, as a sequence of $(n-1)$ index pairs $(m_i,p_i)$, or, alternatively, the values $(\pi_{m_i}, \pi'_{p_i})$; $1\le m_1<m_2<\ldots<m_{n-1}\le n$, $1\le p_1<p_2<\ldots<p_{n-1}\le n$, for which $\pi'_{p_i}-\pi_{m_i}\in\{0,1\}$.  There exist $(n-1)$ {\it necessarily non-intersecting} straight lines that connect the $(\pi_{m_i}, \pi'_{p_i})$ pairs. These lines are either vertical or oblique, and the shape of the entire LCM looks either like
$$\vert\enspace\vert\ldots\vert\enspace/\enspace/\ldots/\enspace\vert\enspace\vert\enspace\ldots\vert,$$ or like$$\vert\enspace\vert\ldots\vert\enspace\backslash\enspace\backslash\ldots\backslash\enspace\vert\enspace\vert\ldots\vert.$$  The mismatched positions in $\pi,\pi'$, i.e., those that do {\it not} belong to an $(m_i,p_i)$ pair, may be (i) such that they can be connected without intersecting any other line (this occurs when there are no oblique lines); or (ii) immediately to different sides of a matched pair (there is one oblique line), or (iii) farther apart (at least two oblique lines), as seen in the examples below.  In the first example,

$$\pi:\quad3\quad8\quad7\quad4\quad2\quad9\quad6\quad10\quad5\quad1$$
$$\pi':\quad3\quad4\quad8\quad5\quad2\quad9\quad7\quad10\quad6\quad1,$$
\centerline{\it Example (i)}

\medskip

\noindent  the mismatch could be considered to be between the 8 (in $\pi$) and the 4 (in $\pi'$); or we have
$$\pi:\quad2\quad1\quad9\quad3\quad4\quad7\quad8\quad6\quad5$$
$$\pi':\quad2\quad1\quad9\quad3\quad7\quad4\quad8\quad6\quad5,$$

\centerline{\it Example (ii)}

\medskip

 \noindent where the mismatch could be viewed as being between the two 7s or between the two 4s; and finally in the last example 

$$
\pi:\quad1\quad2\quad6\quad9\quad8\quad3\quad7\quad4\quad5$$
$$\pi':\quad1\quad2\quad7\quad9\quad3\quad8\quad6\quad4\quad5,$$

\centerline{\it Example (iii)}

\noindent the 8 in $\pi$ and the 6 in $\pi'$ are mismatched.  Example (iii) contains the same permutations as the ones introduced earlier; utilizing the notion of LCMs and mismatches, we have the $(\pi_{m_i}, \pi'_{p_i})$ pairs consisting of $$(1,1), (2,2), (6,7), (9,9), (3,3), (7,8), (4,4), (5,5).$$
The skipped symbols, 8 in $\pi$ and 6 in $\pi'$, provide the mismatch -- and the corresponding covering $(n+1)$-permutation can be constructed from the above matching by introducing the $\emptyset$ symbols to appear below or above the mismatched numbers exactly as before. Thus, the two ways of representing co-coverable permutations -- using either the $\emptyset$ symbols or LCMs, are equivalent, and we will use one or the other as appropriate.

Given a mismatch in which the two numbers are the same, say $a$, they may be represented by either $a$ or $a+1$ in the covering $(n+1)$-permutation.  This gives two coverings depending on where $a$ is placed.  But the common mismatched number may be chosen in up to two ways, as in Example (ii), for a total of 4 covering permutations; this maximal possibility exists if and only if $\pi'$ is obtained from $\pi$ by a single swap of adjacent elements.  
If the two mismatched numbers $a,b$ are different, the case where these are ``atop" each other (as in Example (i)) yields two covering permutations, since the $\emptyset$ symbols may be inserted in the two alternate forms
$${{a}\atop{\emptyset}}\enspace{{\emptyset}\atop{b}}\quad{\rm or}\quad  {{\emptyset}\atop{b}}\enspace{{a}\atop{\emptyset}},$$
which lead to different covering $\rho$s.  If the mismatched numbers are farther apart and different, then it is clear that there is only one covering permutation.

In what follows, we will denote the positions of potential mismatches by $(\alpha,\beta), (\gamma,\delta)$, and set
$$a=\pi_\alpha, b=\pi'_\beta; c=\pi_\gamma, d=\pi'_\delta.$$  We now need to show that all possibilities have been covered, in other words, that there do not exist two LCMs with mismatched pairs at positions $(\alpha,\beta)$ and at $(\gamma,\delta)$, where $c,d$ are not obtained from $a,b$ by the processes described in the previous paragraph, and where $c,d$ yield genuinely different covering permutations.  First observe that alternative choices of $c$ and $d$ {\it do} exist.  For instance, in Example (i), $c$ and $d$ may be taken to be 7 and 4; or 7 and 3; or 8 and 3, but these {\it do not} yield $\rho$s  that are different from the ${{\emptyset}\atop{4}}\enspace{{8}\atop{\emptyset}}$ alignment.

We proceed by contradiction, and by considering all possible ways in which the mismatched pairs $(a,b)$ and $(c,d)$ can interact.  First note that $a, b$ must satisfy $a\ge b$ in order the counteract the fact that all $\pi'_{p_i}\ge \pi_{m_i}$  for all $i$.  Moreover, if $a=b$ then the LCM matches identical elements.  

\medskip

\noindent CASE 1:  Mismatches $a$ and $b$ are atop each other; mismatches $c$ and $d$ are also atop each other; $c$ and $d$ are to the right of $a$ and $b$ respectively.  

Since $a$ and $b$ are atop each other, we must have $\pi^{-1}(i)=\pi'^{-1}(i)$ for each $i\le b-1; i\ge a+1$; $\pi^{-1}(i)\ne\pi'^{-1}(i)$ for each $b\le i\le a$ (assuming without loss of generality that $a\ne b$).  Similarly, $i$ and $i$ are matched for each $i\le d-1; i\ge c+1$ and not for other $i$s.  This forces $[b,a]\cap[d,c]=\emptyset$, and thus $\pi=\pi'$, a contradiction.

\medskip

\noindent CASE 2:  $a$ and $b$ are atop each other; $c$ and $d$ are to the right of $a$ and $b$ respectively.

Assume that $d$ is to the right of $c$.  $c$ and $d$ are initially matched to elements below and above it respectively.  When the mismatched pair becomes $(c,d)$, all matches between $c$ and $d$ are tilted to the right, e.g., the number below the $c$ is now matched to the number to the right of $c$. Since $c\ge d$, and since $i\in \pi'$ may only be matched with either $i$ or $i-1$ in $\pi$, we see that the sequences of integers between $\gamma$ and $\delta$ are identical in $\pi$ and $\pi'$, and consist of a monotone decreasing block of consecutive integers. It follows that all integers lower than $d$ in $\pi$ are matched to themselves in $\pi'$, as are all integers higher than $c$.  This forces both $\pi$ and $\pi'$ to be the same permutation of the identity, a contradiction.

\medskip

\noindent CASE 3:  $\alpha,\beta,\gamma,\delta$ are in increasing order from left to right.

We proceed as in Case 2, except that there are now two monotone blocks --  between $a$ and $b$; and between $c$ and $d$.

\medskip 

\noindent CASE 4: The positions appear in the order $\alpha,\gamma,\beta,\delta$ or one of its variants.

Assume that the $(a,b)$ LCM has oblique lines that go from bottom left to top right.  We see that the oblique lines corresponding to the $(c,d)$-mismatch LCM start between $\alpha$ and $\beta$ but continue beyond $\beta$.   If, however, $\gamma\ge\alpha+1$, or $\delta\ge\beta+1$, we see that there are again one or two monotone blocks, forcing a contradiction.  Thus $\alpha=\gamma$ and $\beta=\delta$, which reduces us to one mismatch.  We thus need to consider when $\alpha=\beta=\gamma$ or $\alpha=\beta=\delta$, as in Example (i), with the (7,4) or (8,3) mismatches.  In this case, either $\alpha$ has to be compatible with both $\beta,\delta$; or $\beta$ has to be compatible with both $\alpha,\gamma$.  The same covering permutations thus result with both sets of mismatches.

\medskip

\noindent CASE 5: The positions appear in the order $\gamma,\alpha,\beta,\delta$ or one of its variants.

Here a contradiction arises if $\gamma\le\alpha-2$ or $\delta\ge\beta+2$.  Assume, therefore, that $\gamma=\alpha-1; \delta=\beta+1$; $\alpha$ and $\beta$ could be arbitrarily far apart, but the point is that the $a-b$ mismatch could be rewired to become a $c-b$ or $a-d$ or even a $c-d$ mismatch, as illustrated in Example (i), with the 7-3 mismatch.  The fact that $\alpha$ and $\beta$ could be arbitrarily far apart is illustrated by the example below:

$$\pi:\quad5\quad6\quad9\quad10\quad11\quad12\quad4\quad7\quad8\quad1\quad2\quad3$$
$$\pi':\quad6\quad9\quad10\quad11\quad12\quad3\quad5\quad7\quad8\quad1\quad2\quad4,$$
\centerline{\it Example (iv)}

\noindent where the (6,3) mismatch can be rewired to become a (5,3) mismatch without changing the covering permutation.  If $\gamma=\alpha-1; \delta=\beta+1$; a $c-d$ rewiring may lead to up to two entirely new $\rho$s, as in Example (ii).

\medskip

\noindent This completes the proof of Lemma 5.\end{proof}

%the lowest $\ell$ elements $1,2,\ldots, \ell$ and the highest $n-h+1$ elements $h,h+1,\ldots,n$ enjoy a perfect match.  Keep in mind that the $\emptyset$ symbols have not yet been placed.  Also, at least one of $\ell, n-h+1$ is strictly positive since, if not, the covering $(n+1)$-permutation would not contain either a 1 or an $(n+1)$.  Next, there must be a certain {\it non-negative} number of positions in which $\sigma'$ is (coordinate-wise) one larger than $\sigma$  -- and then a mismatched situation, i.e., indices $i$ for which $\pi'_i-\pi_i\notin\{0,1\}$.  It is impossible for the cardinality of such $i$s to be more than twoThere are two cases to consider.

We are now ready to prove the main result of this section.

\begin{thm} Consider the probability model in which each $\pi\in S_{n+1}$ is independently picked with probability $p$.  Let the resulting random collection of permutations be denoted by $\ca$.  Then,
\[p\le\frac{\log n}{n}(1+o^*(1))\Rightarrow \p(\ca\ {\rm is\ a\ cover\ of\ }\S_n)\to0\enspace(n\to\infty),\]
and
\[p\ge\frac{\log n}{n}(1+o(1))\Rightarrow \p(\ca\ {\rm is\ a\ cover\ of\ }\S_n)\to1\enspace(n\to\infty).\]
\end{thm}
\begin{proof}We use the Janson inequalities, see  \cite{as}.  By the lower Janson inequality, we have, with $X$ denoting as before the number of uncovered $\pi\in S_n$,
\begin{eqnarray*}\p(X=0)&\ge&\prod_{i=1}^{n!}\p({\rm permutation}\ i\ {\rm is\ covered})\\
&\ge&\prod_{i=1}^{n!}1-(1-p)^{n^2}\\
&\ge&\prod_{i=1}^{n!}\exp\{-(1-p)^{n^2}/(1-(1-p)^{n^2})\}\\
&=&\exp\{-n!(1-p)^{n^2}/(1-(1-p)^{n^2})\to1\\
\end{eqnarray*}
if $\e(X)\to0$, or, using Stirling's formula, if $C\sqrt{n}\lr\frac{n}{e}\rr^ne^{-n^2p}\to0$, i.e., if $$p=\frac{\log n}{n}-\frac{1}{n}+\frac{\log n}{2n^2}+\frac{\omega(n)}{n^2}=\frac{\log n}{n}(1+o(1)),$$
where $\omega(n)\to\infty$ is arbitrary.  This proves the second part of the result.  To prove the first part, we invoke the upper Janson inequality to give
\begin{equation}
\p(X=0)\le\exp\{- \l+\Delta\}
\end{equation}
where $\l=\e(X)$ and
\[\Delta=\sum_i\sum_{j\sim i}\p(i\ {\rm and\ }j\ {\rm are\ both\ uncovered}),\]
and $i\sim j$ if permutations $i$ and $j$ can be covered by the same $(n+1)$-permutation.  Since, by Lemma 5, for any $i$ there are at most $n^3$ permutations $j$ that can be jointly covered with $i$, and, if this is the case, there are at most $4$ covering $(n+1)$-permutations, it follows that 
\[\Delta=n!n^3(1-p)^{2n^2+2-4},\]
so that (2) yields
\[\p(X=0)\le\exp\{-n!(1-p)^{n^2+1}+n!n^3(1-p)^{2n^2-2}\}\to0\]
if $n!(1-p)^{n^2}\to\infty$, i.e., if
\begin{equation}p=\frac{\log n-1+\frac{1}{2}\frac{\log n}{n}-\frac{\omega(n)}{n}}{n}=\frac{\log n}{n}(1+o^*(1)),\end{equation}
where $\omega(n)\to\infty$ is arbitrary.  This completes the proof of Theorem 6.
\end{proof}
\section{Poisson Approximation in the ``Gap"}  The proof of Theorem 6 reveals that $\e(X)$ undergoes a rapid transition when $p$ is around the level given by (3).  In fact, if  
\[p=\frac{\log n-1+\frac{1}{2}\frac{\log n}{n}-\frac{K}{n}}{n}, K\in{\mathbb R},\]
then for large $n$, $\e(X)\sim{\sqrt{2\pi}}e^{-K}$ and $\p(X=0)\sim\exp\{-{\sqrt{2\pi}}e^{-K}\}$.  Much more is true, however, as we shall show next:  The entire probability distribution $\cl(X)$ of $X$ can be approximated, in the total variation sense, by that of a Poisson random variable with mean $\l=\e(X)$ in a range of $p$s that allows for large means.  This result can be thought of as being a probabilistic counterpart to Theorem 6, and is proved using the Stein-Chen method of Poisson approximation \cite{bhj}:
\begin{thm} Consider the model in which each $\pi\in\S_{n+1}$ is independently chosen with probability $p$, thus creating a random ensemble $\ca$ of $(n+1)$-permutations.  Then
$\tv(\cl(X),\P(\l))\to0$ if $p\ge \frac{\log n}{n^2}(1+\epsilon)$, where $\P(\l)$ denotes the Poisson distribution with parameter $\l$, $\epsilon>0$ is arbitrary, and 
the total variation distance $\tv$ is defined by
\[\tv(\cl(Y),\cl(Z))=\sup_{A\subseteq {\mathbb Z^+}}\vert\p(Y\in A)-\p(Z\in A)\vert.\]
\end{thm}
\begin{proof}As before, we set $X=\sum_{j=1}^{n!} I_j$, with $I_j=1$ if the $j$th $n$-permutation is uncovered by permutations in $\ca$ ($I_j=0$ otherwise), and let $\l=\e(X)=n!(1-p)^{n^2+1}$.  Consider the following coupling, for each $j$:  If $I_j=1$, we ``do nothing," setting $J_i=J_{ji}=I_i, 1\le i\le n!$.  If, on the other hand, the $j$th permutation is covered by one or more $(n+1)$-permutations in $\ca$, we ``deselect" these permutations, setting $J_i=1$ if the $i$th permutation is uncovered after this change is made; $J_i=0$ otherwise.  Now it is clear that $J_i\ge I_i$ for $i\ne j$, since there is no way that an uncovered permutation can magically get covered after a few $(n+1)$-permutations are deselected.  Also, setting $N=n!$, we have for each $j$,
\[\cl(J_1,J_2,\ldots, J_N)=\cl(I_1,\ldots,I_N\vert I_j=1).\]
Corollary 2.C.4 in \cite{bhj} thus applies, telling us that 
\begin{equation}\tv(\cl(X),\P(\l))\le\frac{1-e^{-\l}}{\l}\lr\v(X)-\l+2\sum_j\p^2(I_j=1)\rr,\end{equation}
Bounding $1-e^{-\l}$ by one, (4) yields
\begin{equation}\tv(\cl(X),\P(\l))\le\frac{\v(X)}{\l}-1+2(1-p)^{n^2+1}.\end{equation}
The last term in (5), namely $(1-p)^{n^2+1}$ can easily be verified to tend to zero as long as $p\gg 1/n^2$, so we turn to a computation of $\v(X)$:
\[\v(X)=\sum_j\lc\e(I_j)-\e^2(I_j)\rc+\sum_{i\sim j}\{\e(I_iI_j)-\e(I_i)\e(I_j)\},\]
so that 
\begin{eqnarray*}\frac{\v(X)}{\l}-1&\le&\frac{\sum_{i\sim j}\{\e(I_iI_j)-\e(I_i)\e(I_j)\}}{\l}\\
&\le&n^3\lr(1-p)^{n^2-3}-(1-p)^{n^2+1}\rr\\
&=&n^3(1-p)^{n^2+1}\{(1-p)^{-4}-1\}\\
&\le&5pn^3\exp\{-n^2p\}\to0
\end{eqnarray*}
provided that $p\ge A\frac{\log n}{n^2}$, with $A>1$.\end{proof}
\section{Open Problems}  Proving the best possible result along the lines of Conjecture 4 would clearly be a most desirable outcome of this work.  Can, e.g. we write ``$\lim$" instead of ``$\lim\sup$"?  Is $A=1$?  Secondly, can the first bound in Lemma 5 be improved?  (Note, however, that such an improvement would provide only marginal improvements in Theorems 6 and 7.)  Finally, extending the results of this paper to encompass coverings of $n$ permutations by $(n+k)$-permutations  would lead to several interesting questions and new techniques; for example, Theorem 2 can readily be generalized for $k\ge 2$.  Also, the maximum degree of the dependency graph induced by the indicator random variables $\{I_i\}$ is low even for $k\ge2$, so that these variables are almost independent, and one may even envision tight results as $k\to \infty$.   
\section{Acknowledgment}  The research of all four authors was supported by NSF Grant 1004624.


\begin{thebibliography}{99}
\bibitem{as} N.~Alon and J.~Spencer, 1992.  {\it The Probabilistic Method,} John Wiley, New York.
\bibitem{bhj} A.~Barbour, L.~Holst, and S.~Janson, 1992.  {\it Poisson Approximation,} Oxford University Press.
\bibitem{es} P.~Erd\H os and J.~Spencer, 1974.  {\it Probabilistic Methods in Combinatorics,} Academic Press, New York
\bibitem{fs} P.~Flajolet and R.~Sedgewick, 2009.  {\it Analytic Combinatorics,} Cambridge University Press.
\bibitem{gj} A.~Godbole and S.~Janson (1996). ``Random covering designs," {\it J. Combinatorial Theory, Series A} {\bf 75}, 85--98. 
\bibitem{gs} A.~Godbole and P.~Sissokho (1999). ``Near-matches and successions in random permutations," {\it Congressus Numerantium} {\bf 135}, 159--170. 
\bibitem{gss} A.~Godbole, D.~Skipper, and R.~Sunley (1996). ``$t$-covering arrays: upper bounds and Poisson approximations," {\it Combinatorics, Probability and Computing} {\bf 5}, 105--118.  
\bibitem{gtv} A.~Godbole, S.~Thompson, and E.~Vigoda (1996). ``General upper bounds for covering numbers," {\it Ars Combinatoria} {\bf 42}, 211--221.
\bibitem{mw} A.~Myers and H.~Wilf (2003). ``Some new aspects of the coupon collector's problem," {\it SIAM J.~Discrete Mathematics} {\bf 17}, 1--17.
\bibitem{ps} N.~Pippenger and J.~Spencer (1989).  ``Asymptotic behaviour of the chromatic index for hypergraphs," {\it J.
Combinatorial Theory, Series A} {\bf 51}, 24--42.
\bibitem{ro} V.~R\"odl (1985).  ``On a packing and covering problem," {\it European J. Combinatorics} {\bf 5}, 
69--78.
\bibitem{sp} J.~Spencer (1995). ``Packing and covering via a branching process," {\it Rand. Structures Algorithms} {\bf 7}, 167--172.
\end{thebibliography}
\end{document}